\documentclass[11pt,a4]{article}
\input{amssym}
\title{{\Large \bf The Annihilating-Ideal  Graph  of Commutative  Rings II\footnote{ The research
 of the first author was in part supported by
a grant from IPM (No. 87160026).}}}
\author{{\large\bf M. Behboodi$^{a,b}$\footnote{Corresponding
author.} and Z. Rakeei$^a$ }
\\{\footnotesize $^{a}$Department of Mathematical Science,
 Isfahan University of Technology,  } \\{\footnotesize P. O. Box: 84156-83111, Isfahan, Iran }\\{\footnotesize $^{b}$School of Mathematics,
  Institute for Research in Fundamental Sciences (IPM), } \\{\footnotesize P. O. Box: 19395-5746, Tehran, Iran}\\{\footnotesize
mbehbood@cc.iut.ac.ir, \ \ sanamzhr@yahoo.com}}

\def\be{\begin{enumerate}}
\def\ee{\end{enumerate}}

\newtheorem{ttheo}{Theorem}[section]
\newtheorem{ccoro}[ttheo]{Corollary}
\newtheorem{llem}[ttheo]{Lemma}
\newtheorem{eexam}[ttheo]{Example}
\newtheorem{rrem}[ttheo]{Remark}
\newtheorem{ppro}[ttheo]{Proposition}

\newenvironment{pproof}{\noindent{\bf Proof. }}{}

\date{}
\input{epsf}
 \begin{document}
  \maketitle
\begin{center}
  {\small\bf Abstract}
  \end{center}

  {\small In this paper we continue our study of  annihilating-ideal graph of commutative rings, that was introduced in Part I (see [5]).
  Let $R$ be  a commutative ring with ${\Bbb{A}}(R)$ its set of ideals  with nonzero annihilator and $Z(R)$ its set of zero divisors.
 The annihilating-ideal graph of $R$ is defined  as the (undirected) graph ${\Bbb{AG}}(R)$ that its vertices are
  $\Bbb{A}(R)^* =\Bbb{A}(R)\hspace{-1mm}\setminus\{(0)\}$ in which for every  distinct  vertices  $I$ and $J$,
   $I\hspace{-0.6mm}-\hspace{-1.7mm}-\hspace{-1.7mm}-\hspace{-0.5mm}J$ is an edge if and only if
  $IJ=(0)$.   First, we study  the diameter of ${\Bbb{AG}}(R)$. A complete  characterization for the possible diameter is given
   exclusively in terms of the ideals of $R$ when either $R$ is a
   Noetherian ring  or $Z(R)$ is not an ideal of $R$.  Next, we study coloring of  annihilating-ideal graphs.
 Among other results, we characterize when either $\chi({\Bbb{AG}}(R))\leq 2$ or $R$ is reduced and $\chi({\Bbb{AG}}(R))\leq \infty$.
 Also it is shown that for each  reduced ring $R$,  $\chi(\Bbb{AG}(R))= cl(\Bbb{AG}(R))$. Moreover, if
 $\chi(\Bbb{AG}(R))$ is finite, then $R$ has a finite number of
 minimal primes, and if $n$ is this number, then $\chi(\Bbb{AG}(R))= cl(\Bbb{AG}(R))= n$. Finally, we show that for a Noetherian
  ring $R$,  $cl(\Bbb{AG}(R))$ is finite if and only if
for every ideal $I$ of $R$ with $I^2=(0)$, $I$ has finite number
of $R$-submodules.}  \vspace{2mm}
   \vspace{4mm}\\

 \noindent{\bf 0. Introduction}
 \vspace{3mm}

The present paper is a sequel to [5]  and so the notations
introduced in Introduction  of [5] will remain in force. Thus
throughout the paper, $R$ denotes a  commutative ring with
identity, $Z(R)$ denotes the the set of all zero divisors of $R$
and $\Bbb{I}(R)$ denotes the set of all proper ideals of $R$. If
$X$ is either an element or a subset of $R$, then the annihilator
of $X$ is $Ann(X)=\{r\in R | rX=0\}$. The zero divisor graph of
$R$, denoted by $\Gamma(R)$, is a graph with the vertex set
$Z(R)^*:=Z(R)$ in which for every two vertices $x$ and $y$,
$x\hspace{-0.5mm}-\hspace{-1.5mm}-\hspace{-1.5mm}-\hspace{-0.5mm}y$
 is an edge if and only if $x\neq y$ and $xy=0$. As [5], we say that the ideal $I$ of $R$ is  an
{\it annihilating-ideal} if $Ann(I)\neq (0)$ (i.e., there exists a
nonzero ideal $J$ of $R$ such that  $IJ=(0)$). Let  $\Bbb{A}(R)$
be the set of all annihilating-ideals of $R$. Then the {\it
annihilating-ideal graph} of  $R$, denoted by ${\Bbb{AG}}(R)$, is
a undirected simple  graph   with the vertex set
$\Bbb{A}(R)^*:=\Bbb{A}(R)\setminus\{(0)\}$ in which every two
distinct   vertices  $I$ and $J$ are adjacent if and only if
$IJ=(0)$ (see Part I [5] for more details).

   Recall that a graph $G$ is connected if there is a
path between every  two distinct vertices. For distinct vertices
$x$ and $y$ of $G$, let $d(x, y)$ be the length of the shortest
path from $x$ to $y$ and if there is no such path we define
$d(x,y)=\infty$. The  diameter of $G$ is $diam(G) = sup \{d(x,y) \
: \ x \ and \ y \ are \ distinct \ vertices \ of \ G\}$. The girth
of $G$, denoted by $gr(G)$, is defined as the length of the
shortest cycle in $G$ and $gr(G)=\infty$ if $G$ contains no
cycles. A graph in which each pair of distinct vertices is joined
by an edge is called a complete graph. Also, if a  graph $G$
contains one vertex to which all other vertices are joined and $G$
has no other edges, is called a star graph. In [3, Theorem 2.1],
it is shown that for every  ring $R$, ${\Bbb{AG}}(R)$ is a
connected graph and $diam(\Bbb{AG})(R)\leq 3$,  and if
${\Bbb{AG}}(R)$ contains a cycle, then $gr({\Bbb{AG}}(R))\leq 4$
(see [7]). In Section 1 of this paper, we study the diameter of
the annihilating-ideal graphs. By using the papers  [5] and [6],
we determine the relationship between the diameter of
${\Bbb{AG}}(R)$ and $\Gamma(R)$. In particular, a complete
characterization for the possible diameter is given exclusively in
terms of the ideals of $R$ when either $R$
 has finitely many minimal primes  or $Z(R)$ is not an ideal.

 A clique of a graph is a maximal complete subgraph and the number of vertices
in the largest clique of graph $G$, denoted by $cl(G)$, is called
the clique number of $G$. Let $\chi(G)$ denote the chromatic
number of the graph $G$, that is, the minimal number of colors
needed to color the vertices of $G$ so that no two adjacent
vertices have the same color. Obviously $\chi(G) \geq cl(G)$. Beck
in [4] conjectured that $\chi(\Gamma(R))=cl(\Gamma(R))$. But, D.D.
Anderson and M. Naseer gave a counterexample to this conjecture in
[2]. In fact, the  counterexample is the ring $R=\Bbb{Z}_4[X, Y,
Z]/(X^2-2, Y^2-2, Z^2, 2X, 2Y, 2Z, XY, XZ, YZ-2)$ for which
$\chi(\Gamma(R))=5$ but $cl(\Gamma(R))=4$ (note that in Beck's
coloring all elements of the ring were vertices of the graph but
the vertices of $\Gamma(R)$ are nonzero zero divisors of $R$). In
Section 2, we look at the coloring of the annihilating-ideal graph
of rings. First, we show that for the counterexample above we have
$\chi(\Bbb{AG}(R))=cl(\Bbb{AG}(R))=4$. On the other hand, however
we have not found any example where
$\chi(\Bbb{AG}(R))>cl(\Bbb{AG}(R))$. The lack of such
counterexamples  together  with the fact that we have been able to
establish the equality  $\chi(\Bbb{AG}(R))=cl(\Bbb{AG}(R))$ for
Anderson-Naseer's counterexample   and also for  reduced  rings
motivates the following conjecture.\vspace{2mm}\\ \noindent{\bf
Conjecture 0.1.} {\it  For every commutative ring $R$,
$\chi(\Bbb{AG}(R))=cl(\Bbb{AG}(R))$.}\vspace{2mm}\\
\indent  In Section 2, among other results, we characterize rings
$R$ for which  $\chi({\Bbb{AG}}(R))\leq 2$. It is shown that for a
reduced  ring  $R$ the following conditions are equivalent: (1)
$\chi({\Bbb{AG}}(R))<\infty$, (2) $cl({\Bbb{AG}}(R))<\infty$, (3)
$\Bbb{AG}(R)$ does not have an infinite clique and (4) $R$ has
finite number of minimal primes. Moreover, if $R$ a non-domain
reduced ring, then $\chi({\Bbb{AG}}(R)$ is the number of minimal
primes of $R$. Also, it is shown that  for a Noetherian   ring
$R$,  $cl(\Bbb{AG}(R)$ is finite if and only if  every ideal $I$
of $R$ with $I^2=(0)$ has finite number of $R$-submodules. Finally
we conjecture that ``if $\Bbb{AG}(R)$ does not have an infinite
clique, then $\chi({\Bbb{AG}}(R))$ is finite".\vspace{2mm}\\

  \noindent{\bf 1. The diameter of an annihilating-ideal graph }\\

   By Anderson and Livingston [3, Theorem 2.3], for every
ring $R$, the zero divisor graph  $\Gamma(R)$ is a connected graph
and $diam(\Gamma(R))\leq  3$. Moreover, if $\Gamma(R)$ contains a
cycle, then gr( $\Gamma(R))\leq 4$ (see [7]). Moreover, Lucas in
[6] characterized the diameter of $\Gamma (R)$ in terms of the
ideals of   $R$. As we have seen in [5,  Theorem 2.1], for a ring
  $R$, $\Bbb{AG}(R)$ is also a connected graph with  $diam(\Bbb{AG}(R)\leq 3$. In general, Lucas's results are not true  for annihilating-ideal
  graphs, but the following proposition more or less summarizes the overall
situation for  the  relationship between the diameter of
  $\Bbb{AG}(R)$ and $\Gamma (R)$.\\

  \noindent{\bf Proposition 1.1.} {\it  Let $R$ be a ring.}\vspace{3mm}\\
  (a) {\it If $diam(\Gamma(R))=0$, then $diam({\Bbb{AG}}(R))=0$.}\\
  (b) {\it If $diam(\Gamma(R))=1$, then $diam({\Bbb{AG}}(R))=0 \ or \ 1$.}\\
  (c) {\it If $diam(\Gamma(R))=2$, then $diam({\Bbb{AG}}(R))=1, 2 \ or \ 3$.}\\
  (d) {\it If $diam(\Gamma(R))=3$, then $diam({\Bbb{AG}}(R))=3$.}\\
(e) {\it If $diam({\Bbb{AG}}(R))=0$, then ${diam(\Gamma(R))=0 \ or \ 1}$.}\\
(f) {\it If $diam({\Bbb{AG}}(R))=1$, then ${diam(\Gamma(R))=1 \ or \ 2}$.}\\
  (g) {\it If $diam({\Bbb{AG}}(R))=2$, then ${diam(\Gamma(R))=2}$.}\\
  (h) {\it If $diam({\Bbb{AG}}(R))=3$, then ${diam(\Gamma(R))=2 \ or \ 3}$.}\\

\noindent{\bf Proof.}  (a). Let   $diam(\Gamma(R))=0$ i.e.,
$\Gamma(R)$ has one vertex. Thus  by [6, Theorem 2.6.],
 $R$ is isomorphic to either $\Bbb{Z}_4$ or $\Bbb{Z}_2[y]/(y^2)$. Clearly, in any case $\Bbb{AG}(R)$ has also one vertex   and so
 $diam({\Bbb{AG}}(R))=0$.\\
  (b). Clearly, if $diam(\Gamma(R))=1$ then $\Gamma(R)$ is a complete graph with more than one vertex.
  Thus  by [6, Theorem 2.6], either (i) $R$ is reduced and isomorphic to $\Bbb{Z}_2\times\Bbb{Z}_2$, or (ii) $R$
  is non-reduced, $Z(R)^2=(0)$ and $R$ is not isomorphic to either $\Bbb{Z}_4$ or
  $\Bbb{Z}_2[y]/(y^2)$. If $R\cong\Bbb{Z}_2\times\Bbb{Z}_2$, then $diam({\Bbb{AG}}(R))=1$. In the later case,   $Z(R)^2=(0)$ implies that
  $\Bbb{AG}(R)$ is also a complete graph (see [5, Theorem 2.7]. Now, if $\Bbb{AG}(R)$ has one vertex
  (as $\Bbb{Z}_{p^2}$ where $p$ is an odd prime number), then
  $diam({\Bbb{AG}}(R))=0$, otherwise, $diam({\Bbb{AG}}(R))=1$.\\
    (c). By [6, Theorem 2.6.], $diam(\Gamma(R))=2$ if and only if either $R$ is reduced with exactly two minimal primes
and at least three nonzero zero divisors, or $Z(R)$ is an ideal
whose square is not (0) and each pair of distinct zero divisors
has a nonzero annihilator. If $diam(\Bbb{AG}(R))=0$, then
${\Bbb{AG}}(R)$ has only one vertex, say $I$, and  $I^2=(0)$. It
follows that $Z(R)=I$ and so by [3, Theorem 2.8], $\Gamma(R)$ is a
complete graph.  Thus $diam(\Gamma(R))=0$ or $1$, contradicting
with our hypothesis. Thus $diam(\Bbb{AG}(R))>0$, i.e.,
$diam(\Bbb{AG}(R))=1,2$ or 3 (since by [5, Theorem 2.1],
$diam(\Bbb{AG}(R))\leq 3$). We claim that  all these three cases
may happen. To see this, let $R_2=F_1\times F_2$  where $F_1$ and
$F_2$ are fields and $F_1\ncong \Bbb{Z}_2$ and $R_1=\Bbb{Z}\times
\Bbb{Z}_2$. Then $diam(\Gamma(R_1))=diam(\Gamma(R_2))=2$ but
$diam({\Bbb{AG}}(R_1))=1$ and $diam({\Bbb{AG}}(R_2))=2$. For an
example of a ring $R$  for which $diam(\Gamma(R))=2$ but
$diam({\Bbb{AG}}(R_2))=3$ see   Example 1.7 of this paper.\\
    (d). Let $diam(\Gamma(R))=3$. Then there exist  $a$, $b\in Z(R)$ such that $d(a,b)=3$. Clearly, $Ra$ and $Rb$ are vertices of
     ${\Bbb{AG}}(R)$ and $RaRb\neq (0)$.     We claim that $Ra\neq Rb$, for if not, then $Ra=Rb$ implies that $a$ and $b$ have the
      same annihilator i.e., $d(a,b)=2$, a contradiction. Now by  [5,  Theorem 2.1], $d(Ra,Rb)\leq 3$ in ${\Bbb{AG}}(R)$. If
$d(Ra,Rb)\neq 3$, then $d(Ra,Rb)=2$, i.e., there exists $0\neq
c\in R$ such that $RcRa=RcRb=(0)$. It follows $d(a,b)=2$ in
$\Gamma(R)$, a contradiction. Thus  $d(Ra,Rb)=3$ and
so $diam({\Bbb{AG}}(R))=3$.\\
 (e). Clearly, $diam({\Bbb{AG}}(R))=0$ if and only if the ring $R$ has only
 one nonzero proper ideal, if and only if $Z(R)$ is the only  nonzero proper
 ideal of $R$. In this case $Z(R)^2=(0)$ and so $\Gamma(R)$ is a complete
 graph. If $\Gamma(R)$ has one vertex, then by [6, Theorem 2.6 (1)], $R$ is isomorphic to either $\Bbb{Z}_4$ or
 $\Bbb{Z}_2[y]/(y^2)$ and so  $diam(\Gamma(R))=0$, otherwise, when $\Gamma(R)$ has more than
 one vertices (as $R=\Bbb{Z}_{p^2}$ where $p$ is an odd prime
 number), then  $diam(\Gamma(R))=1$. \\
 (f). Clearly, $diam({\Bbb{AG}}(R))=1$ if and only if
 $\Bbb{AG}(R)$ is a complete graph with more than one vertex, i.e., either $R=F_1\oplus F_2$, where
$F_1$, $F_2$ are fields, $Z(R)$ is an ideal of $R$ with
$Z(R)^2=(0)$ or $R$ is a local ring with exactly two nonzero
proper ideals $Z(R)$ and $Z(R)^2$ (see [5, Theorem 2.7]). If
$R=\Bbb{Z}_2\oplus \Bbb{Z}_2$ or $Z(R)^2=(0)$, then $\Gamma(R)$ is
a complete graph with more than one vertex i.e.,
$diam(\Gamma(R))=1$. It is easy to check that if $R=F_1\oplus
F_2\not\cong\Bbb{Z}_2\oplus \Bbb{Z}_2$ or $R$ is a local ring with
exactly two nonzero proper ideals $Z(R)$ and $Z(R)^2$, then
$diam(\Gamma(R))=2$.\\
(g). If  $diam({\Bbb{AG}}(R))=2$, then  by (a), (b) and (d) above,
$diam(\Gamma(R))\neq 0, \ 1$ and $3$. Thus $diam(\Gamma(R))=2$.\\
(h). If  $diam({\Bbb{AG}}(R))=3$, then by (a) and (b) above,
$diam(\Gamma(R))\neq 0$ and $1$. Thus  $diam(\Gamma(R))=2$ or $3$.
By (d) above if $diam(\Gamma(R))=3$,  then $diam({\Bbb{AG}}(R))$
is also $3$. For the case $diam({\Bbb{AG}}(R))=3$ but
$diam(\Gamma(R))=2$, see Example 1.7. $\square$\\

In next two theorems, we characterize when $Z(R)$ is not an ideal
of $R$ and either $diam(\Bbb{AG}(R))=2$ or $diam(\Bbb{AG}(R))=3$, respectively. \\

  \noindent{\bf Theorem  1.2.} {\it  Let $R$ be a  ring such that $Z(R)$ is not an ideal of $R$. Then $diam({\Bbb{AG}}(R))=2$ if and only if
  $R$ is reduced with  exactly
two minimal primes and at least three nonzero annihilating ideals.}\\

\noindent{\bf Proof.} $(\Rightarrow)$.   Let
$diam({\Bbb{AG}}(R))=2$. Then $R$ has at least three nonzero
annihilating ideals and  by Proposition 1.1 (g),
${diam(\Gamma(R))=2}$. Now by [6, Theorem 2.6, (3)], $R$ is a
reduced ring  with  exactly two minimal primes. \\
$(\Leftarrow)$. If $R$ is reduced with only two minimal primes
$P_1$ and $P_2$, then each nonzero ideal with a nonzero
annihilator is contained in exactly one of the two minimal primes
and its annihilator is the other minimal prime. This follows
easily from the fact that if $J$ is a nonzero ideal of a reduced
ring, then $Ann(J)$ is a radical ideal equal to the intersection
of the minimal primes that do not contain the $J$. So in this
special case $diam(\Bbb{AG}(R)) = 2$ if and only if at least one
of the minimal primes properly contains a nonzero ideal. Since  $R$ have at least three nonzero annihilating ideals,
$diam(\Bbb{AG}(R)) = 2$.   $\square$\\

\noindent{\bf Theorem 1.3.} {\it Let $R$ be a ring such that
$Z(R)$ is not an ideal. Then the following statements are
equivalents.}\vspace{3mm}\\
(1) {\it $diam({\Bbb{AG}}(R))=3$.}\\
(2)  {\it $diam(\Gamma(R))=3$.}\\
(3) {\it Either $R$ is non-reduced or $R$ is a
reduced ring with more than two min- \indent imal   primes.}\\

\noindent{\bf Proof.} $(1)\Rightarrow(2)$.   Let
$diam({\Bbb{AG}}(R))=3$. If $diam(\Gamma(R))\neq 3$, then by
Proposition 1.1 (h), $diam(\Gamma(R))=2$. Thus by  [6, Theorem
2.6, (3)], $R$ is reduced with  exactly two minimal primes. Since
$diam({\Bbb{AG}}(R))=3$,  there exist at least three nonzero
annihilating ideals in $R$. Thus
 by Theorem 1.2, $diam({\Bbb{AG}}(R))=2$, a contradiction. Therefore,
 $diam(\Gamma(R))=3$.\\
$(2)\Rightarrow(1)$ is by Proposition 1.1, (d).\\
$(2)\Leftrightarrow(3)$ is by [6, Theorem 2.6.]. $\square$\\

Now, we are in position to give a complete  characterization for
the possible diameter of $\Bbb{AG}(R)$ when $Z(R)$ is not an ideal of $R$.\\

\noindent{\bf Theorem 1.4.} {\it Let $R$ be a ring such that
$Z(R)$ is not an ideal of $R$. Then $1\leq diam({\Bbb{AG}}(R)\leq 3$ and} \vspace{3mm}\\
(1) {\it $diam({\Bbb{AG}}(R))=1$ if and only if $R\cong F_1\times
F_2$ where $F_1$ and $F_2$ are fields;}\vspace{1mm}\\
(2) {\it $diam({\Bbb{AG}}(R))=2$ if and only if $R$ is reduced
with  exactly two minimal \indent  primes  and at least three
nonzero
annihilating ideals;}\vspace{1mm}\\
(3) {\it $diam({\Bbb{AG}}(R))=3$ if and only if either $R$ is a
non-reduced  ring or $R$ is \indent a  reduced ring with more than
two minimal   primes.}\\

\noindent{\bf Proof.} Suppose that $Z(R)$ is not an ideal of $R$.
Then by [5, Theorem 2.1], ${diam({\Bbb{AG}}(R))\leq 3}$. If
$diam({\Bbb{AG}}(R)=0$, then $\Bbb{AG}(R)$ has only one vertex,
i.e.,  $Z(R)$ is the only nonzero proper ideal of $R$, a
contradiction. Thus we conclude that $diam({\Bbb{AG}}(R)=1, \ 2$
or $3$.\\
For (1), we note that ${\Bbb{AG}}(R)$ is a complete graph (i.e.,
$diam({\Bbb{AG}}(R))=0$ or  $1$)  if and only if either
$R=F_1\oplus F_2$, where $F_1$, $F_2$ are fields, or $Z(R)$ is an
ideal of $R$, $Z(R)^3=(0)$ and for each ideal $I\subsetneqq Z(R)$,
$IZ(R)=(0)$ (see [5, Theorem 2.7]). Thus $diam({\Bbb{AG}}(R))=1$
if and only if $R\cong F_1\times F_2$ where $F_1$ and $F_2$ are fields.\\
Finally statements  (2) and (3) are from Theorems 1.2 and 1.3
above. $\square$\\

In the next theorem,  we provide a sufficient condition for
$\Bbb{AG}(R)$  to have diameter $3$ when $R$  is a non-reduced
ring. First we need the following lemma (see also [6, Lemma 2.3]).\\

\noindent{\bf Lemma 1.5.} {\it Let $R$ be a ring and let $I$ be an
annihilating-ideal  of $R$. If $N$ is a nilpotent ideal of $R$,
then  $N+I$ is an annihilating-ideal  of $R$.}\\

\noindent{\bf Proof.} Let $N$ be a nonzero nilpotent ideal and
assume $cI=(0)$  where $c\neq $0. Since $N$ is nilpotent, there is
a positive integer $m$ such that $cN^m=0$ with $cN^{m-1}\neq 0$.
Clearly, $cN^{m-1}\subseteq Ann(N+I)$. $\square$\\

\noindent{\bf Theorem 1.6.} {\it Let $R$ be a non-reduced ring. If
there is a pair of annihilating-ideals $I$, $J$ of $R$  such that
$Ann(I+J)=(0)$, then ${diam({\Bbb{AG}}(R))=3}$.}\\

\noindent{\bf Proof.} Let $I$, $J\in {\Bbb{A}}^*(R)$ be such that
$Ann(I+J)=(0)$. Then $d(I,J)\neq 2$. By Lemma 1.5, neither $I$ nor
$J$ can be nilpotent. If $IJ\neq (0)$, then $d(I,J)=3$. Thus we
may assume $IJ=(0)$. Since $IJ=(0)$, $(I+J)^2=I^2+J^2$ has no
nonzero annihilator. Since $R$ is not reduced, there exists a
nonzero nilpotent $q\in R$. Thus without loss of generality we may
assume that $qJ^2\neq (0)$. Since $I$ is an annihilating-ideal and
$qJ$ is nilpotent, $I+qJ$ is an annihilating-ideal by Lemma 1.5.
On the other hand $I+qJ\neq J$, for if not, then $I\subseteq J$
and so $Ann(I+J)=Ann(J)\neq (0)$, a contradiction. Now consider
the pair $I+qJ$ and $J$. Since $I+J=I+qJ+J$ and $Ann(I+J)=(0)$,
$d(I+qJ, J)\neq 2$. But $(I+qJ)J=qJ^2\neq (0)$. Thus $d(I+qJ,
J)=3$ and
$diam({\Bbb{AG}}(R))=3. \ \square$\\

\noindent{\bf Example  1.7.} Let $R$ be the ring in [6, Example
5.5.]. Then $R$ is a non-reduced ring  and there is a pair of
annihilating-ideals $I$, $J$ of $R$  such that $I+J$ is not an
annihilating-ideal (see [6, Example 5.5, (6)]. Thus  by Theorem
1.6,  $diam({\Bbb{AG}}(R))=3$, but since $R$ is a McCoy ring (see
[6, Example 5.5, (4)] $diam(\Gamma(R))=2$. (Note: a ring $R$ is
said to be a McCoy ring if each finitely generated ideal contained
in $Z(R)$ has a nonzero annihilator).\\

It is easy to see that  a reduced ring does have the property that
an ideal with a nonzero annihilator must be contained in at least
one minimal prime.\\

\noindent{\bf Lemma  1.8.} {\it Let $R$ be a reduced   ring with
finitely many minimal primes. If $R$ has  more than two minimal
primes,  then $diam({\Bbb{AG}}(R))=3$.}\\

\noindent{\bf Proof.} Since  $R$  has more than  two minimal
primes, it is not a domain.  Let $I$ be a nonzero  ideal of $R$
with nonzero annihilator. Since $R$ is a reduced ring, $I+Ann(I)$
is not an  annihilating-ideal. Let $J=Ann(I)$. Since   $IJ=(0)$,
each minimal prime of $R$ contains at least one of $I$ and $J$,
but not both of them (since $R$ has finitely many minimal primes,
each minimal primes has nonzero annihilator).  Thus without loss
of generality we may assume there are minimal primes $P_1$, $P_2$
and $P_3$ such that $I\subseteq P_1\cap P_2$, $I\not\subseteq
P_3$, $J\subseteq P_3$, $J\not\subseteq P_1$ and $J\not\subseteq
P_2$. Let $q\in P_1\cap P_3\setminus P_2$ and consider the pair
$I+qJ$ and $J$. Since $R$ is reduced, $IJ=(0)$ and neither $I$ nor
$q$ is contained in $P_2$, $(0)\neq qJ^2=J(I+qJ)$. Since
$I+qJ\subseteq P_1$,  that $I+qJ$ is a nonzero annihilating-ideal
of $R$. On the other hand,
$I+J=I+qJ+J$ is an ideal with no nonzero annihilator. Thus $d(I+qJ, J)=3$ and $diam({\Bbb{AG}}(R))=3$. $\square$\\

Now, we are in position to give a complete  characterization for
the possible diameter of $\Bbb{AG}(R)$ when $R$ has  finitely many
minimal primes (as Noetherian rings).\\

\noindent{\bf Theorem  1.9.} {\it Let $R$ be a ring  with finitely
many minimal primes. Then }\vspace{3mm}\\
(1) {\it  $diam({\Bbb{AG}}(R))=0$ if and only if $Z(R)$ is the only nonzero proper ideal \indent of $R$.}\vspace{1mm}\\
(2) {\it $diam({\Bbb{AG}}(R))=1$ if and only if  $R=F_1\oplus
F_2$, where $F_1$, $F_2$ are  fields, \indent $Z(R)$ is an ideal
of $R$ with  $Z(R)^2=(0)$ or $R$ is a local ring with
  exactly  \indent two nonzero proper ideals
$Z(R)$ and $Z(R)^2$.}\vspace{1mm}\\
 (3) {\it $diam({\Bbb{AG}}(R))=2$
if and only if either $R$ is reduced with exactly  two
 \indent minimal primes and at least three nonzero
annihilating-ideals, or $R$ is  \indent not reduced, $Z(R)$ is an
ideal whose square is not $(0)$ and  for  each pair \indent  of
annihilating-ideals  $I$  and  $J$, $I+J$ is an annihilating-ideal.}\vspace{1mm}\\
(4) {\it  $diam({\Bbb{AG}}(R))=3$ if and only if either  $R$ is
reduced with  more
 than two \indent  minimal primes   or $R$ is non-reduced and there are annihilating-ideals \indent  $I\neq J$,
such  that $I+J$ is not an
annihilating-ideal and either.}\\

\noindent{\bf Proof.} (1) is evident. \\
The statement (2) is from [5, Theorem 2.7].\\
(3) $(\Rightarrow)$. Suppose that $diam({\Bbb{AG}}(R))=2$. By
Proposition 1.1 $diam(\Gamma(R))=2$. Thus by [6, Theorem 2.6 (3)],
either (i) $R$ is reduced with exactly two minimal primes and at
least three nonzero zero divisors, or (ii) $Z(R)$ is an ideal
whose square is not (0) and each pair of distinct zero divisors
has a nonzero annihilator. The Case (i) implies that $Z(R)$ is not
an ideal of $R$ and so  by Theorem 1.2, $R$ is a reduced ring with
exactly two  minimal primes and at least three nonzero
annihilating-ideals. Assume that the Case (ii) holds. If $R$ is
reduced  with  more than two  minimal primes, then  by Lemma 1.8,
$diam({\Bbb{AG}}(R))=3$, a contradiction. Thus we can assume that
$R$ is not reduced, and so by Theorem 1.6, for each pair of
annihilating-ideals $I$ and $J$, $I+J$ is an annihilating-ideal. \\
(3) $(\Leftarrow)$ is evident.\\
 Finally statement (4) is from Theorems 1.3,  1.6, Lemma 1.8 and
 statement  (3) above. $\square$\\

\noindent{\bf Remark 1.10.} If $R$ is a reduced ring with more
than two minimal ideals, then by [1, Theorems 2.2 and 2.4],
$gr(\Gamma(R))=3$. This fact is also true for the
annihilating-ideal graph of a reduced ring $R$ with finitely many
minimal primes $P_1, P_2,\cdots, P_n$ ($n\geq 3$). In fact, if
$P_1, P_2,\cdots, P_n$ are  distinct minimal primes,  then
$P_2P_3\cdots
P_n-\hspace{-2mm}-\hspace{-2mm}-\hspace{-2mm}-P_1P_3\cdots
P_n-\hspace{-2mm}-\hspace{-2mm}-\hspace{-2mm}-P_1P_2\widehat{P_3}\cdots
  P_n-\hspace{-2mm}-\hspace{-2mm}-\hspace{-2mm}-P_2P_3\cdots P_n$ is a cycle of length
  $3$.   We have not found any
examples of a ring reduced ring $R$  with infinite minimal primes
 for which $gr(\Bbb{AG}(R))\neq 3$. The lack of such
counterexamples, together  above fact, motivates the following
 conjecture:\\

\noindent{\bf Conjecture 1.11}. {\it Let $R$ be  a reduced ring
with more than two minimal primes. Then $gr(\Bbb{AG}(R))=3$.}\\

\noindent{\bf 2. Coloring of the  annihilating ideal graphs }\\

The goal of this section is to study of coloring of the
annihilating ideal graphs of rings, in particular  the interplay
between $\chi(\Bbb{AG}(R))$ and $cl(\Bbb{AG}(R))$.\\

Beck [4] showed that if $R$ is a reduced ring or a principal ideal
ring, then $\chi(\Gamma(R))=cl(\Gamma(R))$ He also showed that for
$n=1, 2$ or $3$, $\chi(\Gamma(R))=n$ if and only if
$cl(\Gamma(R))=n$. Based on these positive results, Beck
conjectured that $\chi(\Gamma(R))=cl(\Gamma(R))$ for each ring $R$
when $\chi(\Gamma(R))<\infty$. In [2], D.D. Anderson and M. Naseer
gave a counterexample to this conjecture. We start this section
with the following interesting result about  Anderson-Naseer's counterexample.\\

\noindent{\bf Proposition  2.1.} {\it Let
\begin{center}
$R=\Bbb{Z}_4[X, Y, Z]/(X^2-2, Y^2-2, Z^2, 2X, 2Y, 2Z, XY, XZ,
YZ-2).$\end{center}
  Then $\chi(\Bbb{AG}(R))=cl(\Bbb{AG}(R))=4$.}\\

\noindent{\bf Proof.} Clearly, $R$ is a finite local ring with 32
elements and  $J(R)=Z(R)=\{0, 2, x+2, y, y+2, x+y, x+y+2, z, z+2,
x+z, x+z+2, y+z, y+z+2, x+y+z, x+y+z+2\}$. One can easily see that
 the  ideals
$(2)$, $(x)$, $(y)$, $(z)$, $(x+y)$, $(x+z)$, $(y+z)$, $(x+y+z)$,
$(x, y)$,  $(x, z)$, $(y, z)$,
$(x, y+z)$, $(y, x+z)$, $(z, x+y)$, $(x, y, z)$  are the only  nonzero proper ideals of $R$.\\
Thus the graph  $\Bbb{AG}(R)$ has 15 vertices and so  is the graph
in Figure 1
below.\\\\\\\\\\\\\\\\\\\\\\\\\\\\\\\\\\\\\\\\\\\\\\\\

Now by using  Figure 1,  it easy to see that $\{(2), (x), (y),
(y+z)\}$ is a maximal clique and each other clique of
$\Bbb{AG}(R)$ has at most 3 elements. Thus $cl(\Bbb{AG}(R))=4$.
Also, the Figure 1,  shows that the minimal number of colors
needed to color $\Bbb{AG}(R)$ is  4 colors, which
we label as 1, 2, 3 and 4. Thus $\chi(\Bbb{AG}(R))=cl(\Bbb{AG}(R))=4$. $\square$\\

Next, we characterize rings $R$  for which $\chi (R)=1$ or 2.\\

\noindent{\bf Proposition  2.2.} {\it Let $R$ be a ring. Then
$\chi(\Bbb{AG}(R))=1$ if and only if $R$ has
only one nonzero proper  ideal.}\\

 \noindent{\bf Proof.} Since $\Bbb{AG}(R)$ is a  connected  graph and
 $\chi(\Bbb{AG}(R))=1$, it can not have more than one vertex, and so by [5, Theorem 1.4], $R$ has
only one nonzero proper  ideal. The converse is clear.   $\square$\\

\noindent{\bf Theorem 2.3.} {\it For a  ring $R$  the
following statements are equivalent:}\vspace{3mm}\\
 (1) {\it $\chi(\Bbb{AG}(R))=2$.}\\
 (2) {\it $\Bbb{AG}(R)$ is a  bipartite graph with two nonempty parts.}\\
 (3) {\it $\Bbb{AG}(R)$ is a complete bipartite graph with two nonempty parts.}\\
 (4) {\it Either  $R$ is a reduced ring with exactly two
 minimal primes or   $\Bbb{AG}(R)$ is \indent a  star graph with more than one vertex.}\\

 \noindent{\bf Proof.} $(1)\Leftrightarrow(2)$ and $(3)\Rightarrow(2)$ are clear.\\
 $(2)\Rightarrow(4)$. Suppose that $\Bbb{AG}(R)$ is bipartite with
two parts $V_1$ and $V_2$.\\
 \noindent{Case 1}: The ring $R$ is a
reduced ring.  Suppose $R$ has at least three minimal prime ideals
and let $P_1$, $P_2$ and $P_3$ be three distinct minimal primes.
Let $b_1\in P_1\cap P_2\setminus P_3$. Since both $R_{P_1}$ and
$R_{P_2}$ are fields, neither $P_1$ nor $P_2$ contains $Ann(b_1)$.
So by Prime Avoidance, there is an element $c\in Ann(b_1)\setminus
P_1 \cup P_2$. Next let $d\in P_1\setminus P_2\cup P_3$ and set
$b_2 = cd\in P_1\cap P_3\setminus P_2$. Since $c\in Ann(b_1)$, the
ideals $b_1R$ and $b_2R$ form an edge. The ideal $b_1R+b_2R$ is a
finitely generated ideal and contained in the minimal prime $P_1$,
thus it has a nonzero annihilator $b_3$. But then both $\{b_1R,
b_3R\}$ and $\{b_2R, b_3R\}$ are edges. Hence $\Bbb{AG}(R)$ is not
bipartite. Thus  $R$ is a reduced ring with exactly two
 minimal primes. \\
\noindent{Case 2}: The ring  $R$ is not reduced. Assume $x^2=(0)$,
where $0\neq x\in R$. Without loss of generality we can assume
that $Rx\in V_1$.  We claim that either $Rx$ is a minimal ideal of
$R$ or for each  $0\neq z\in R$ such that $Rz\subsetneq Rx$, $Rz$
is a minimal ideal of $R$. To see this let $(0)\neq Rz_1\subsetneq
Rx$,  $(0)\neq Rz_1\subset Rx$,$(0)\neq Rz_2\subset Rz_1$ for some
$z_1, z_2\in R$. Since $(Rx)^2=(0)$, $Rz_1Rx=Rz_2Rx=(0)$, and so
$Rz_1$, $Rz_2\in V_2$ since $\Bbb{AG}(R)$ is bipartite. But
$Rz_1Rz_2\subseteq (Rx)^2=(0)$, a contradiction. Thus without loss
of generality we can assume that $Rx$ is a minimal ideal of $R$.
Thus $P=Ann(x)$ is a maximal ideal of $R$ and also $Rx\subseteq
P$. We claim that $V_1=\{Rx\}$. Since $x^2 = 0$, every prime
contains $xR$. So if $P = xR$, then it is the only nonzero ideal
of $R$ (under the assumption $xR$ is a minimal ideal), and hence
$\Bbb{AG}(R)$ is not  a  bipartite graph, a contradiction. Thus
$P\neq Rx$. Since $PRx=(0)$ and $\Bbb{AG}(R)$ is bipartite, $P\in
 V_2$. Now let $V_1\neq \{Rx\}$ and $L\in V_1\setminus \{Rx\}$. If $LP=(0)$, $LRx\subseteq LP=(0)$, a
 contradiction. Thus $LP\neq (0)$ and since $\Bbb{AG}(R)$ is connected, there
exists $K\in V_2$ such that $K\neq P$ and $LK=(0)$. It follows
that $L\cap P\neq (0)$, $(L\cap P)Rx=(0)$. Therefore $L\cap P\in
V_2$. But $(L\cap P)K=(0)$, a contradiction. Thus in each case
$V_1=\{Rx\}$ and it follows that $\Bbb{AG}(R)$ is a star graph with more than one vertex.\\
 $(4)\Rightarrow(3)$. Clearly if $\Bbb{AG}(R)$ is a star graph, then $\Bbb{AG}(R)$ is a complete bipartite graph.
 Thus we assume that $R$ is a reduced ring with  exactly two  minimal primes,  say  $P_1$
and $P_2$. Then $P_1\cap P_2=(0)$ and  $Z(R)=P_1\cup P_2$. Thus if
$I$ is a nonzero ideal with a nonzero annihilator, by  Prime
Avoidance Theorem (see [8, Theorem 3.61], $I$ is contained in one
of $P_1$ and $P_2$ and its annihilator is the other prime.
Therefore,   if $I_1,I_2\subseteq P_i$ ($i=1, 2$), then
$I_1I_2\neq (0) $, and hence $\Bbb{AG}(R)$ is bipartite. On the
other hand, for vertices $I$, $J$ such that $I\subseteq P_1$ and
$J\subseteq P_2$, we have $IJ\subseteq P_1\cap P_2=(0)$.
Therefore, $\Bbb{AG}(R)$ is a complete bipartite graph with two nonempty parts. $\square$\\

\noindent{\bf Corollary  2.4.} {\it Let $R$ be an Artinian ring.
Then the following statements are equivalent:}\vspace{3mm}\\
 (1) {\it $\chi(\Bbb{AG}(R))=2$.}\\
 (2) {\it $\Bbb{AG}(R)$ is a  bipartite graph with two nonempty parts.}\\
 (3) {\it $\Bbb{AG}(R)$ is a complete bipartite graph with two nonempty parts.}\\
 (4) {\it Either $R\cong F_1\times F_2$ for some fields  $F_1$ and $F_2$ or $R$ is a local ring such \indent  that
  $\Bbb{AG}(R)$ is a star graph with more than one vertex.}\vspace{2mm}\\

 \noindent {\bf Proof.} By Theorem 2.3, $(1)\Leftrightarrow
 (2)\Leftrightarrow (3)$.\\
 $(2)\Rightarrow (4)$. Assume that $\Bbb{AG}(R)$ is a  bipartite
 graph with two nonempty parts. By Theorem 2.3, either  $R$ is a reduced ring with exactly two
 minimal primes  or  $\Bbb{AG}(R)$ is a  star graph with more than one vertex. In the first
 case  $R\cong F_1\times F_2$ for some fields  $F_1$ and $F_2$.
 In the second case, if $R$ is a local ring, we are done. Otherwise,  then by [5, Theorem 2.6],   $R\cong F_1\times F_2$ for some fields
  $F_1$ and $F_2$. \\
$(4)\Rightarrow (2)$ is clear. $\square$\\

In [5, Corollary 2.4], it is shown that for a reduced ring $R$,
$\Bbb{AG}(R)$ is a star graph if and only if $R\cong F\oplus D$,
where $F$ is a field and $D$ is an integral domain. It follows
that $R$ has exactly two
 minimal primes $F\times (0)$ and $(0)\times D$. Also, for a reduced ring $R$,  $\chi(\Gamma(R))=2$ if and only
 if $R$ has  exactly two  minimal primes (see [4, Theorem 3.8]).  Thus by using these  facts  and  Theorem 2.3,  we have
the following interesting result for reduced rings.\\

 \noindent{\bf Corollary 2.5.} {\it Let $R$ be a reduced ring. Then
the following statements are equivalent:}\vspace{3mm}\\
 (1) {\it $\chi(\Bbb{AG}(R))=2$.}\\
 (2) {\it $\chi(\Gamma(R))=2$.}\\
 (3) {\it $\Bbb{AG}(R)$ is a  bipartite graph with two nonempty parts.}\\
 (4) {\it $\Bbb{AG}(R)$ is a complete bipartite graph with two nonempty parts.}\\
 (5) {\it $R$ has  exactly two
 minimal primes.}\vspace{2mm}\\

One can easily see that for every ring $R$, $$\chi(\Gamma(R))=2 \
\  \Longrightarrow \ \chi(\Bbb{AG}(R))=2.$$
 But the following example
shows that in general, if $R$ is not reduced,
$$\chi(\Bbb{AG}(R))=2 \ \not\Longrightarrow  \ \chi(\Gamma(R))=2.$$

\noindent{\bf Example 2.6.}  Let $R=\Bbb{Z}_{p^3}$ where $p$ is a
positive  prime number. Then it is easy to check that
$\chi(\Bbb{AG}(R))=2$, but $\chi(\Gamma (R))= p$ (see [4,
Proposition 2.3] in which, it is known that both the chromatic
number and the clique number are $p+1$, since
Beck includes $0$ as a vertex).\\

The following proposition shows that if $\chi(\Gamma(R))$ is
finite, then $\chi(\Bbb{AG}(R))$ is also finite.\\

\noindent{\bf Proposition  2.7.} {\it Let $R$ be a   ring. If
$\Bbb{AG}(R)$ has an infinite clique, then $\Gamma(R)$ has an
infinite clique.}\\

\noindent {\bf Proof.} Suppose that $C$ is an infinite clique of
$\Bbb{AG}(R)$. Let $\Upsilon$ be the set of all vertices $I$ of
$C$ such that $I^2=(0)$.\\
\noindent{Case 1}: If $\bigcup \Upsilon$ is infinite, then for
every $x,y\in
\Upsilon$, $xy=(0)$ and so $\Gamma(R)$ has an infinite clique.\\
\noindent{Case 2}: If $\bigcup \Upsilon$ is finite. then there are
finitely many vertices $I$ of $C$ such that $I^2=(0)$. Thus
$C\setminus \Upsilon$ is infinite. We claim that in this case, for
every $I\in C\setminus\Upsilon$, $I\nsubseteq \bigcup_{\phi\in
\Phi}J_{\phi}$, where $\{J_{\phi} | \phi\in \Phi\}$ is an
arbitrary subset of $C\setminus \Upsilon$ such that
$I\notin\{J_{\phi} | \phi\in \Phi\}$. Because if $I\subseteq
\bigcup_{\phi\in \Phi}J_{\phi}$, then $I^2=(0)$, a contradiction.
Now, let $\{I_i | i\in \Bbb{N}\}$ be a subset of distinct vertices
of $C\setminus\Upsilon$. Let $x_1\in I_1$. Since $I_2\nsubseteq
I_1$, there exists $x_2\in I_2\setminus I_1$ such that $x_2\neq
x_1$. Again, since $I_3\nsubseteq (I_1\cup I_2)$, there exists
$x_3\in I_3\setminus (I_1\cup I_2)$ such that $x_3\neq x_1,x_2$.
Continuing
this way, we can get an infinite clique of $\Gamma(R)$. $\square$\\

We are going to characterize reduced rings $R$ for which
$\chi(\Bbb{AG}(R))$ is finite. To see that, we need some
prefaces.\\

\noindent {\bf Lemma 2.8.} {\it  Let $R$ be a reduced ring such
that $\Bbb{AG}(R)$ does not have an infinite clique. Then $R$ has ACC on ideals of the form $Ann(I)$ where $I$ is an ideal of $R$.}\\

\noindent{\bf Proof.} Let $Ann(I_1)\subset Ann(I_2)\subset
Ann(I_3)\subset ...$ be a chain in $\Bbb{A}(R)$.  Clearly,
$I_iAnn(I_{i+1})\neq (0)$ for each $i\geq 1$. Thus for each $i\geq
1$, there exists $x_i\in I_i$ such that $x_i Ann(I_{i+1})\neq (0)
$. Let $J_i=x_i Ann(I_{i+1}), i=2, 3, ...$. Then if $i\neq j$,
$J_i\neq J_j$
because if $J_i= J_j$, then $J_i^2= J_j^2= 0$, contradiction. $\square$\\

\noindent {\bf Lemma 2.9.} ( [4, Lemma 3.6.]) {\it  Let $P_1 =
Ann(x_1)$ and $P_2 = Ann(x_2)$ be two distinct elements of
$Spec(R)$. Then we have $x_1x_2= 0$.}\\

\noindent {\bf Theorem 2.10.} {\it For a reduced ring $R$, the
following statements are equivalent:}\\
1) {\it $\chi(\Bbb{AG}(R))$ is finite.}\\
2) {\it $cl(\Bbb{AG}(R))$ is finite.}\\
3) {\it $\Bbb{AG}(R)$ does not have an infinite clique.}\\
4) {\it $R$ has finite number of minimal prime ideals.}\\

\noindent{\bf Proof.} $(1)\Rightarrow (2)\Rightarrow (3)$ is clear.\\
  $(4)\Rightarrow (1)$. Let $(0) = P_1\cap P_2\cap
...\cap P_k$, where $P_1, ..., P_k$ are prime ideals. Define a
coloring  $f$ on $\Bbb{A}(R)^*$ by putting $f(J)= min \{n\in
\Bbb{N}:
J\nsubseteq P_n\}$ for $J\in \Bbb{A}(R)^*$. Then $\chi(\Bbb{AG}(R))\leq k$.\\
$(3)\Rightarrow (4)$. Suppose that $R$ doesn't have an infinite
clique. So by Lemma 2.8, $R$ has ACC on ideals of the form
$Ann(I)$,  $I\in{\Bbb{I}}(R)$. Thus the set $\{Ann(x): 0\neq x\in
R\}$ has maximal ideals, and it is easy to see that these are
prime ideals of $R$. Let  $Ann(x_\lambda)$, $\lambda\in\Lambda$ be
the different maximal members of the family $\{Ann(x): 0\neq x\in
R\}$.
 By Lemma 2.9, the index set $\Lambda$ is  finite.
 Pick $x\in R$, $x\neq 0$. Then $Ann(x)\subseteq Ann(x_\lambda)$ for some $\lambda\in \Lambda$.
Now $x_\lambda x\neq 0$ because otherwise, $x_{\lambda}^2= 0$,
which is a contradiction. Thus  $x\not\in Ann(x_\lambda)$. This
means
$\bigcap_{\lambda\in \Lambda} Ann(x_\lambda)= (0). \ \square$\\

Now we are in position to prove Conjecture 0.1 in the case of
reduced rings. In fact the next corollary is an immediate
conclusion of Theorem 2.10.\\

\noindent {\bf Corollary 2.11.} {\it Let $R$ be a reduced ring.
Then $\chi(\Bbb{AG}(R))= cl(\Bbb{AG}(R))$. Moreover, if
$\chi(\Bbb{AG}(R))$ is finite, then $R$ has a finite number of
minimal primes, and if $n$ is this number, then
$$\chi(\Bbb{AG}(R))= cl(\Bbb{AG}(R))= n.$$}

Adding Corollary 2.11, and [4, Theorem 3.9.] to  Theorem 2.10, the next corollary is indisputable.\\

\noindent {\bf Corollary 2.12.} {\it Let $R$ be a nonzero reduced
ring. Then $\chi(\Bbb{AG}(R))$ is finite if and only if
$\chi(\Gamma (R))$ is finite. Furthermore, if they are finite,
then $$\chi(\Bbb{AG}(R))=\chi(\Gamma (R)).$$}

Now, we will have a characterization of Noetherian rings $R$ for
which  $cl(\Bbb{AG}(R))<\infty $. Before that, we need  the
following evident
 lemma.\\

\noindent {\bf Lemma 2.13.} {\it Let $R$ be a   ring. If
$cl(\Bbb{AG}(R))$ is finite, then for every nonzero ideal $I$ of
$R$ with $I^2=(0)$, $I$
has finite number of $R$-submodules.}\\

\noindent {\bf Theorem 2.14.} {\it Let $R$ be a Noetherian
  ring. Then $cl(\Bbb{AG}(R))$ is finite if and only if
for every ideal $I$ of $R$ with $I^2=(0)$, $I$ has finite number
of $R$-submodules.}\\

\noindent {\bf Proof.} ($\Rightarrow $) is by Lemma 2.13.\\
($\Leftarrow$ ) Suppose that for every ideal $I$ of $R$ with
$I^2=(0)$, $I$ has finite number of $R$-submodules. Let $C$ be a
largest clique in $\Bbb{AG}(R)$, and let $\Upsilon$ be the set of
all vertices $I$ of $C$ with $I^2=(0)$.  If
$\Upsilon\neq\emptyset$, then $J=\sum_{I\in \Upsilon} I$ is again
a vertex of $C$ and $J^2=(0)$. So by our hypothesis, $J$ has
finite number of $R$-submodules. But if $I\in \Upsilon$, every
$R$-submodule of $I$ is an $R$-submodule of $J$. Thus for every
$I\in \Upsilon$, $I$ have finite number of $R$-submodules and
specially $\Upsilon$ has finite element. We claim that $C\setminus
\Upsilon$ has finite member, too. Suppose that $\{I_1, I_2, ...\}$
is an infinite subset of $C\setminus \Upsilon$. Consider the chain
$I_1\subseteq I_1+I_2\subseteq I_1+I_2+I_3\subseteq ...$. Since
$R$ is Noetherian, there exists $n\in \Bbb{N}$, such that $I_1+
...+I_n=I_1+ ...+I_n+I_{n+1}$, i.e. $I_{n+1}\subseteq I_1+
...+I_n$. Since by our choice of $C$, $I_1+ ...+I_n\in C$,
$I_{n+1}$ is adjacent to $I_1+ ...+I_n$ and thus $I_{n+1}^2=(0)$,
a contradiction. Thus $C$ has finite number
of vertices and from there, $cl(\Bbb{AG}(R))$ is finite. $\square$\\

By Theorem 2.10, one may naturally guess the following statements
for the general case.We close this section with this conjecture.\\

\noindent {\bf Conjecture 2.15.} {\it Let $R$ be a ring. If
$\Bbb{AG}(R)$ does not have an infinite clique, then
$cl(\Bbb{AG}(R))$ is finite.}\\

\noindent{\bf  Acknowledgments}\vspace{2mm}

   \noindent{ This work was partially supported by the Center of Excellence of Algebraic
Methods and Application of Isfahan University of Technology. }\vspace{2mm}\\

\noindent{\bf References}\vspace{2mm}\\
{\footnotesize
 \noindent$[1]$ D. D. Anderson, S. B.Mulay,  On the diameter and girth of a zero-divisor graph, J. \indent  Pure Appl. Algebra 210 (2007)
 543-550).\\
 $[2]$ D. D. Anderson, M. Naseer,  Beck's coloring of a commutative
  ring, J. Algebra 159 \indent (1993)  500-514.\\
$[3]$ D. F. Anderson, P. S. Livingston,  The zero-divisor graph of
a  commutative  ring, J. \indent Algebra  217 (1999) 434-447.\\
$[4]$ I. Beck,  Coloring of commutative   rings, J. Algebra 116
(1988) 208-226. \\
$[5]$ M. Behboodi, Z. Rakeei, The annihilating-ideal  graph  of commutative  rings I, submit- \indent ted.\\
 $[6]$ T. G. Lucas,  The diameter of a zero divisor graph. J.
Algebra 301 (2006) 174-193.\\
 $[7]$ S. B. Mulay,  Cycles and symmetries of zero-divisors, Comm.
Algebra, 30 (2002) 3533- \indent 3558.\\
  $[8]$  R. Y. Sharp, Steps in   algebra, Second edition, London Mathematical Society Student \indent Texts,  51.
  Cambridge University Press, Cambridge, 1990.
  \end{document}